\documentclass[12pt]{amsart}
\usepackage{latexsym}
\usepackage{amsmath, amssymb}
\usepackage{palatino}
\usepackage{color}
\usepackage[all]{xy}

\newtheorem*{thm*}{Theorem}

\begin{document}

\title{The Jiang-Su algebra does not always embed}
\author{Marius Dadarlat}
\address{Department of Mathematics, Purdue University,
150 N. University St., West Lafayette, \indent IN, 47907-2067, USA}
\email{mdd@math.purdue.edu}
\author{Ilan Hirshberg}
\address{Department of Mathematics, Ben Gurion University of the Negev, P.O.B. 653, Be'er \indent Sheva 84105, Israel}
\email{ilan@math.bgu.ac.il}
\author{Andrew S. Toms}
\address{Department of Mathematics and Statistics, York University,
4700 Keele St., Toronto, \indent Ontario, Canada, M3J 1P3}
\email{atoms@mathstat.yorku.ca}
\author{Wilhelm Winter}
\address{School of Mathematical Sciences, University of Nottingham, University Park, \indent Nottingham, NG7 2RD, United Kingdom}
\email{wilhelm.winter@nottingham.ac.uk}

\date{\today}

\thanks{The authors were partially supported by the Fields Institute; M.D. was partially supported by NSF grant \#DMS-0500693;
I.H. was partially supported by the Israel Science Foundation
(grant No. 1471/07); A.T. was partially supported by NSERC}

\begin{abstract}
We exhibit a unital simple nuclear non-type-I C$^*$-algebra into which the
Jiang-Su algebra does not embed unitally.  This answers a question of M. R{\o}rdam.
\end{abstract}

\maketitle


The Jiang-Su algebra, denoted by $\mathcal{Z}$ (\cite{js}), occupies a central position in
the structure theory of separable amenable C$^*$-algebras.  The
property of absorbing the Jiang-Su algebra tensorially is a necessary, and, in
considerable generality, sufficient condition for the confirmation of G. A. Elliott's
$\mathrm{K}$-theoretic rigidity conjecture for simple separable amenable C$^*$-algebras (\cite{r},
\cite{wl}).
The uniqueness question for this algebra is therefore of great interest.  M. R{\o}rdam
observed that if $\mathcal{C}$ is a class of unital separable C$^*$-algebras, and
$A \in \mathcal{C}$ has the properties that (i) for every $B \in \mathcal{C}$ there is
a unital $*$-homomorphism $\gamma:A \to B$ and (ii) every unital $*$-endomorphism of
$A$ is approximately inner, then $A$ is the only such algebra, up to isomorphism.
(This follows from an application of Elliott's Intertwining Argument.)  Every unital $*$-endomorphism
of $\mathcal{Z}$ is approximately inner (\cite{js}), and there are no obvious obstructions to the existence
of a unital $*$-homomorphism $\gamma:\mathcal{Z} \to A$ for any unital separable C$^*$-algebra
$A$ without finite-dimensional quotients.  Indeed, such a $\gamma$ always exists when $A$ has real rank zero,
and examples show that the existence of $\gamma$ is strictly weaker than tensorial absorption of
$\mathcal{Z}$--see \cite{er} and \cite{t}, respectively.  All of this begs the question, first posed by R{\o}rdam:
"Does every unital C$^*$-algebra without finite-dimensional quotients admit a unital embedding
of $\mathcal{Z}$?", see \cite{er}.  We prove that the answer is negative, even
when the target algebra is simple and nuclear.

\begin{thm*}
 There is a unital simple nuclear infinite dimensional
 C$^*$-algebra (in fact, an AH algebra)
into which the Jiang-Su algebra does not embed unitally.
\end{thm*}
In the remainder of the paper we give some
background discussion and prove the theorem.
For a pair of relatively prime integers $p,q>1$, we set
\[
Z_{p,q} = \{ f \in \mathrm{C}([0,1]; \mathrm{M}_{p} \otimes \mathrm{M}_q) \ | \
f(0) \in \mathrm{M}_p \otimes 1, f(1) \in 1 \otimes \mathrm{M}_q \}.
\]
Each $Z_{p,q}$ is contained unitally in $\mathcal{Z}$.  If a unital C$^*$-algebra $A$
admits no unital $*$-homomorphism $\gamma:Z_{p,q} \to A$, then there is no unital embedding
of $\mathcal{Z}$ into $A$.

Let $A$, $B$ be unital C$^*$-algebras, and let $e,f \in A$ be
projections satisfying $(n+1)[e] \leq n[f]$ in the Murray-Von Neumann semigroup $\mathrm{V}(A)$ for
some $n \in \mathbb{N}$.  It is implicitly shown  in the proof of
\cite[Lemma 4.3]{rordam stable rank} that if
$\gamma:Z_{n,n+1} \to B$ is a unital $*$-homomorphism,
then $[e\otimes 1_B] \leq [f\otimes 1_B]$ in $\mathrm{V}(A\otimes B)$; tensor products are minimal.  Example 4.8 of \cite{hrw} exhibits a sequence
$(B_j)_{j \in \mathbb{N}}$ of unital separable C$^*$-algebras with
the following property:  there are projections $e,f \in B_1 \otimes B_2$
such that $4[e]\leq 3[f]$, but $[e \otimes 1_{\bigotimes_{j=3}^nB_j}]
\not \leq [f \otimes 1_{\bigotimes_{j=3}^nB_j}]$ for any $n \geq
3$.  Using R{\o}rdam's result, one concludes that there is no
unital $*$-homomorphism $\gamma:Z_{3,4} \to \bigotimes_{j=3}^n B_j$
for any $j \geq 3$.  (In fact, there is nothing special about $Z_{3,4}$.
A similar construction can be carried out for a wide variety of $Z_{p,q}$s.)

To simplify notation, we renumber the $B_j$s so that
$Z_{3,4}$ does not embed into $\bigotimes_{j=1}^{n}B_j$
for any $n \in \mathbb{N}$.
  For each $i \in \mathbb{N}$, set $D_i =
\bigotimes_{j=1}^{i} B_j$. We will perturb the canonical embeddings
\[
\psi_i := \mathrm{id} \otimes 1_{B_{i+1}}: D_i \longrightarrow D_i \otimes B_{i+1} = D_{i+1}
\]
to maps $\phi_i$ with the property that $(D_{i},\phi_i)_{i \in
\mathbb{N}}$ has simple limit $D$.  Any such limit, simple or not,
fails to admit a unital $*$-homomorphism $\gamma:Z_{3,4} \to D$,
and so also fails to admit a unital embedding of $\mathcal{Z}$.
Indeed, suppose that such a $\gamma$ did exist.  Then, by the
semiprojectivity of $Z_{3,4}$ (\cite{js}), there would exist a
unital $*$-homomorphism $\tilde{\gamma}: Z_{3,4} \to D_i$ for some
$i$, contradicting our choice of $D_i$.  We remark that, in
particular, $\bigotimes_{j=1}^{\infty} B_j$ admits no unital
embedding of $\mathcal{Z}$. This algebra is a continuous field of
C$^*$-algebras whose fibres are $\mathcal{Z}$-absorbing -- in
fact, its fibres are all isomorphic to the CAR algebra (see
\cite[Example 4.8]{hrw}).

The $B_j$s have the form $(e_j \oplus f_j)\left(\mathrm{C}(X_j) \otimes \mathcal{K} \right)(e_j \oplus f_j)$,
where $e_j$ and $f_j$ are rank one projections and $X_j=(\mathrm{S}^2)^{\times m(j)}$.
Let $\alpha: X_j\to X_j$
be a homeomorphism homotopic to the identity map, and view $B_j$ as a corner of $\mathrm{C}(X_j) \otimes \mathrm{M}_n$
for some sufficiently large $n \in \mathbb{N}$.  The map $\alpha$ induces an automorphism $\alpha^*$ of
$\mathrm{C}(X_j) \otimes \mathrm{M}_n$, $\alpha^*(f)= f\circ \alpha$.  In
general, $\alpha^*$ will not carry $B_j$ into $B_j$, but this can be corrected.  Since $\alpha$ is homotopic
to the identity, the projection $e_j \oplus f_j$ is homotopic, and hence unitarily equivalent, to its image
under $\alpha^*$.  If $u$ is a unitary implementing this equivalence, then $\overline{\alpha} := (\mathrm{Ad}(u)
\circ \alpha^*)|_{B_j}$ is an automorphism of $B_j$.  For our purposes, the salient property of $\overline{\alpha}$
is this:  if $f \in B_j$ and $f(x) \neq 0$ for some $x \in X_j$, then $\overline{\alpha}(f)(\alpha^{-1}(x)) \neq 0$.

It remains to construct the $\phi_i$, and prove the simplicity of the resulting inductive limit algebra $D$.
Let us set $Y_i := \prod_{j=1}^{i} X_j$ where $i\in \mathbb{N}$ or $i=\infty$. We endow $Y_i$  with the metric
$d(x,y)=\sum_{j=1}^i 2^{-j}d_j(x^j,y^j)$ where $d_j$ is the canonical metric on $X_j$,  normalized so that $X_j$ has diameter equal to one.

 Choose a dense sequence $(z_k)_{k \in \mathbb{N}}$ in
$Y_\infty$.  Fix $z_0 \in Y_\infty$ and for each $k\in \mathbb{N}$ let $\beta_k:Y_\infty\to Y_\infty$ be a cartesian product
of isometries  of $X_j$s which are homotopic to the identity and such that $\beta_k(z_0)=z_k$.
Let $(\alpha_k)_{k\in \mathbb{N}}$ be an enumeration of the set $\{\beta_m\beta_n^{-1}:n,m\in \mathbb{N}\}$. It is easy to see that
for any point $x\in Y_\infty$ and any $i\in \mathbb{N}$, the sequence
$(\alpha_k(x))_{k\geq i}$ is dense in $Y_\infty$.
Note that each $\alpha_k$ is also a  cartesian product
of isometries $\alpha_k^j$ of $X_j$ homotopic to $\mathrm{id}_{X_j}$. Let us set $\alpha_{k,[i]}=\prod_{j=1}^i \alpha_k^j$.
Let $\pi_i:Y_\infty\to Y_i$ be the co-ordinate projection. Then $\pi_i\alpha_k(x)=\alpha_{k,[i]}(\pi_i(x))$ for $x\in Y_\infty$.
 Therefore
 for any point $y\in Y_i$, the sequence
$(\alpha_{k,[i]}(y))_{k\geq i}$ is dense in $Y_i$. So is the sequence $((\alpha_{k,[i]})^{-1}(y))_{k\geq i}$ since
   each $\alpha_{k,[i]}$ is an isometry.  By the compactness of $Y_i$  it follows   that
for any nonempty open set $U$ of $Y_i$, there is $j\geq i$ such that
$Y_i=\bigcup_{k=i}^j (\alpha_{k,[i]})^{-1}(U)$.

For each $i\leq k \in \mathbb{N}$, let $\overline{\alpha_{k,[i]}}:D_i \to D_i$ be the automorphism induced, in the manner described
above, by the homeomorphism $\alpha_{k,[i]}:Y_i\to Y_i$.

Observe that the canonical embedding $\psi_i:D_i \to D_{i+1}$ is the direct sum of two non-unital embeddings:
\[
\psi_i^{(1)} \stackrel{\mathrm{def}}{=} \mathrm{id} \otimes e_{i+1}: D_i \to D_i \otimes e_{i+1} \subseteq D_{i+1},
\]
and
\[
\psi_i^{(2)} \stackrel{\mathrm{def}}{=} \mathrm{id} \otimes f_{i+1}: D_i \to D_i \otimes f_{i+1} \subseteq D_{i+1}.
\]
Set $\phi_i^{(1)} = \psi_i^{(1)}$, and
\[
\phi_i^{(2)} \stackrel{\mathrm{def}}{=} \overline{\alpha_{i,[i]}} \otimes f_{i+1}: D_i \to D_i \otimes f_{i+1} \subseteq D_{i+1}.
\]
Define $\phi_i:D_i \to D_{i+1}$ to be $\phi_i^{(1)} \oplus \phi_i^{(2)}$.

Let us now verify that $D = \lim_{i \to \infty} (D_i,\phi_i)$ is simple.  It will suffice to prove that for any nonzero
$a \in D_i$ there is some $j \geq i$ such that $\phi_{i,j+1}(a) := \phi_{j} \circ \cdots \circ \phi_i(a)$ is nonzero over every point in the
spectrum of $D_{j+1}$.

For each $\mathbf{v} = (v_i,\ldots,v_{j}) \in \{1,2\}^{j-i+1}$, set
\[
\phi_{i,j+1}^{\mathbf{v}} = \phi_j^{v_{j}} \circ \phi_{j-1}^{v_{j-1}} \circ \cdots \circ \phi_i^{v_i},
\]
and note that $\phi_{i,j+1} = \bigoplus_{\mathbf{v} \in \{1,2\}^{j-i+1}} \phi_{i,j+1}^{\mathbf{v}}$.
For $k \in \{i,\ldots,j\}$, let $\mathbf{v}_k \in \{1,2\}^{j-i+1}$ be the vector which
is equal to 1 in each co-ordinate except the $k^{\mathrm{th}}$ one.    We have (with the exception of the cases $k= i,j$ when the formula reads slightly differently)
\begin{gather*}
\phi_{i,j+1}^{\mathbf{v}_k}(a) = \left[\, \overline{\alpha_{k,[k]}}(a \otimes e_{i+1} \otimes \cdots\otimes e_{k}) \right]
\otimes f_{k+1} \otimes e_{k+2} \otimes \cdots \otimes e_{j+1}
\\ = \left[\, \overline{\alpha_{k,[i]}}(a) \otimes  \overline{\alpha_{k}^{i+1}}(e_{i+1}) \otimes\cdots\otimes \overline{\alpha_{k}^{k}}(e_{k})) \right]
\otimes f_{k+1} \otimes e_{k+2} \otimes \cdots \otimes e_{j+1}.
\end{gather*}
Since $a$ is nonzero on some nonempty open set $U$, the formula above shows that $\phi_{i,j+1}^{\mathbf{v}_k}(a)$ is
nonzero on  $W^k_{i,j+1}:=(\alpha_{k,[i]})^{-1}(U) \times X_{i+1} \times \cdots \times X_{j+1}$,
for any $j\geq i$. As noticed earlier,
there is $j\geq i$ such that
$Y_i=\bigcup_{k=i}^j (\alpha_{k,[i]})^{-1}(U)$.
 Therefore $\bigoplus_{k=i}^{j} \phi_{i,j+1}^{\mathbf{v}_k}(f)$ is nonzero on
$\bigcup_{k=i}^{j} W^k_{i,j+1}=Y_{j+1}=\widehat{D_{j+1}}$,
as required.

\end{document}